\documentstyle[12pt]{article}
\def\Bbb R{{\rm \bf R}}
\def\proclaim#1{\vskip2mm{\bf #1}\em}
\def\endproclaim{\em \vskip2mm}
\def\tag#1{\eqno(#1)}
\def\gathered{\begin{array}{c}}
\def\endgathered{\end{array}}
\def\text{\mbox}

\begin{document}

\title {Reconstruction of a source domain from the Cauchy data}
\author{Masaru IKEHATA\footnote{
Department of Mathematics,
Faculty of Engineering,
Gunma University,
Kiryu 376-8515, JAPAN}}
\date{}
\maketitle

\begin{abstract}
We consider an inverse source problem for the Helmholtz equation in a bounded domain.
The problem is to reconstruct the shape of the support of a source term from the Cauchy data
on the boundary of the solution of the governing equation.
We prove that if the shape is a polygon, one can calculate its support function from such data.
An application to the inverse boundary value problem is also included.

\end{abstract}


\section{Introduction}

In this paper, we consider an inverse source problem, which is more general than traditional
inverse potential problem described in \cite{Is}, initiated by Novikov \cite{N}.

We consider the weak solution $u\in H^1(\Omega)$ of the Helmholtz equation with a source term
$F\in\{H^1_0(\Omega)\}^*$ in a two-dimensional bounded domain $\Omega$ with Lipschitz boundary:
$$\begin{array}{ll}
\displaystyle
\Delta u+k^2u=F
&
\text{in $\Omega$.}
\end{array}
\tag {1.1}
$$
Here $k$ is a fixed non-negative number and known.  We denote by $\nu$ the unit outward normal
to $\partial\Omega$.  Throughout this paper we assume that
$$\displaystyle
\text{supp}\,F\subset\Omega.
$$
Take $\psi\in C^{\infty}_0(\Omega)$ such that $\psi=1$ in a neighbourhood of $\text{supp}\,F$.
Define
$$\begin{array}{ll}
\displaystyle
\tilde{F}(w)=F(\psi\, w)
&
\displaystyle
w\in H^1(\Omega).
\end{array}
$$
This bounded linear functional on $H^1(\Omega)$ does not depend on the choice of such $\psi$.
Using it, we can define the Neumann derivative of $u$ on $\partial\Omega$ as an element of the dual space of
$H^{1/2}(\partial\Omega)$ by the formula
$$\displaystyle
<\frac{\partial u}{\partial\nu}\vert_{\partial\Omega},\,g>
=\int_{\Omega}\nabla u\cdot\nabla w-k^2\int_{\Omega} uw+\tilde{F}(w)
$$
where $g\in H^{1/2}(\partial\Omega)$ and $w$ is a lifting of $g$ in $H^1(\Omega)$.

We call $(u\vert_{\partial\Omega},\,\partial u/\partial\nu\vert_{\partial\Omega})$ the Cauchy data of 
$u$ on $\partial\Omega$.  Roughly speaking, our problem is
to extract some information about the support of $F$ from the Cauchy data of $u$ on $\partial\Omega$.
When $k=0$ and $\Omega$ is a three-dimensional domain, this problem is related to the inverse source
problem in geophysics since the restriction of the gravitational potential to $\Omega$
satisfies the Poisson equation.

The purpose of this paper is to present a simple idea for the reconstruction of the support
of the source in the two-dimensional case because of the simplicity of the situation.
In a subsequent paper we will consider the three-dimensional case.
We consider the special type of the source:
$$\displaystyle
F=\rho(x)\chi_D(x)
\tag {1.2}
$$
or
$$\displaystyle
F=\nabla\cdot\{\rho(x)\chi_D(x)\mbox{\boldmath $a$}\}.
\tag {1.3}
$$
Here $\rho$ is a function defined on an open subset $D$ of $\Omega$ and $\mbox{\boldmath $a$}$ is a non-zero constant vector.
We assume that $D$, $\rho$ and $\mbox{\boldmath $a$}$ are all unknown and Cauchy data $u\vert_{\partial\Omega}$, 
$(\partial/\partial\nu)u\vert_{\partial\Omega}$ are known.
We also assume that 
$$\displaystyle
\overline{D}\subset\Omega
\tag {1.4}
$$
$$\displaystyle
\text{$D$ is a polygon}
\tag {1.5}
$$
$$\displaystyle
\rho\in L^2(D)
\tag {1.6}
$$
and for each vertex $p$ on $\partial D$ there exist an open disk $B(p,\eta)$ centred at $p$ with radius $\eta>0$,
$0<\alpha\le 1$ and a function $\tilde{\rho}\in C^{0,\alpha}(\overline{B(p,\eta)})$ such that
$$\displaystyle
\text{$\rho=\tilde{\rho}$ on $B(p,\eta)\cap\overline{D}$}
\tag {1.7}
$$
and
$$\displaystyle
\text{$\tilde{\rho}(p)\not=0$.}
\tag {1.8}
$$
Recall the support function for $D$:
$$
\displaystyle
h_D(\omega)=\sup_{x\in D}\,x\cdot\omega,\,\omega\in S^1.
$$

{\bf\noindent Definition (regular direction).}
$\omega$ is regular with respect to $D$ if and only if the set
$$\displaystyle
\{x\in\Bbb R^2\,\vert\,x\cdot\omega=h_D(\omega)\}\cap \partial D
$$
consists of only one point.

The main result is as follows.

\proclaim{\noindent Theorem 1.1.}
Assume that $F$ has the form (1.2) or (1.3).  One can calculate
$h_D(\omega)$ for regular $\omega$ with respect to $D$ from the Cauchy data of $u$ on $\partial\Omega$.

\endproclaim

The proof contains how to calculate $h_D(\omega)$ for regular $\omega$.  See (2.12) of Proposition 2.3
in the next section.  Furthermore, since $h_D(\omega)$ is continuous with respect to $\omega$ and the set
of directions which are not regular with respect to $D$ are at most a finite set, from the uniqueness
of $h_D(\omega)$ we automatically obtain two uniqueness theorems about the reconstruction of the
convex hull $[D]$ of $D$.

\proclaim{\noindent Corollary 1.}  Assume that $(D_1,\rho_1)$ and $(D_2,\rho_2)$ are two pairs satisfying (1.4)-(1.8).
Let $F_j$ have the form (1.2) for $(D,\rho)=(D_j,\rho_j)$.  Let $u_j$ be the weak solution of (1.1) for $F=F_j$.
Then
$$\displaystyle
\left(u_1\vert_{\partial\Omega},\frac{\partial u_1}{\partial\nu}\vert_{\partial\Omega}\right)
=\left(u_2\vert_{\partial\Omega},\frac{\partial u_2}{\partial\nu}\vert_{\partial\Omega}\right)
$$
implies
$$\displaystyle
[D_1]=[D_2].
$$

\endproclaim

\proclaim{\noindent Corollary 2.}  Assume that $(D_1,\rho_1,\mbox{\boldmath $a$}_1)$ and 
$(D_2,\rho_2,\mbox{\boldmath $a$}_2)$ are two pairs satisfying (1.4)-(1.8).
Let $F_j$ have the form (1.3) for $(D,\rho,\mbox{\boldmath $a$})=(D_j,\rho_j,\mbox{\boldmath $a$}_j)$.  Let $u_j$ be the weak solution of (1.1) for $F=F_j$.
Then
$$\displaystyle
\left(u_1\vert_{\partial\Omega},\frac{\partial u_1}{\partial\nu}\vert_{\partial\Omega}\right)
=\left(u_2\vert_{\partial\Omega},\frac{\partial u_2}{\partial\nu}\vert_{\partial\Omega}\right)
$$
implies
$$\displaystyle
[D_1]=[D_2].
$$

\endproclaim

Of course, as a direct consequence of Corollaries 1 and 2, we obtain two uniqueness theorems in the inverse potential problem and they seem to be new.
See \cite{Is} for several uniqueness theorems.

Badia and Duong \cite{BD} considered a similar inverse source problem for the Poisson equation
with a product source in a cylindrical geometry.  Using the method of separation of variables,
they find a reconstruction procedure for the shape of the support of the product source under
additional information: the volume or one of the first-order moments of the bottom domain is known.

The method for the proof of Theorem 1.1 is similar to the idea discovered in \cite{E00} and does not require 
a cylindrical structure of the domain.  Therein we considered two variants
of the inverse conductivity problem proposed by Calder\'on \cite{C}.  The first is to reconstruct
the surface of discontinuity of the coefficient $\gamma$ of the equation $\nabla\cdot\gamma\nabla u=0$ in $\Omega$
from the Dirichlet-to-Neumann map; the second is to reconstruct the shape of $\partial D$ on which the Dirichlet data of the solution
of the Helmholtz equation $\Delta u+k^2u=0$ in $\Omega\setminus\overline D$ vanish,
from the Dirichlet-to-Neumann map.  Note that $k$ is fixed.  We found a new application of the so-called complex
plane wave
$$\displaystyle
\exp(x\cdot z),\,z\cdot z=-k^2
$$
to these problems.  Using it and an energy inequality, we reduced the problems to study
the asymptotic behaviour of the integral
$$\begin{array}{ll}
\displaystyle
\int_D\exp\{2\tau(x\cdot\omega-t)\}\,dx
&
\displaystyle
\omega\in S^2
\end{array}
\tag {1.9}
$$
as $\tau\longrightarrow\infty$ for each fixed $t\in\Bbb R$.  We succeeded in calculating the support function for 
the surface from the Dirichlet-to-Neumann map.

The argument in this paper, in contrast, is based on studying the asymptotic behaviour of the integral
$$\begin{array}{ll}
\displaystyle
\int_D\exp\{x\cdot(\tau\omega+i\sqrt{\tau^2+k^2}\,\omega^{\perp})-\tau t)\}\,dx
&
\displaystyle
\omega,\,\omega^{\perp}\in S^1,\,\omega\cdot\omega^{\perp}=0
\end{array}
$$
as $\tau\longrightarrow\infty$ for each fixed $t\in\Bbb R$.  Note that this integral
is a kind of oscillatory integral and has a different character from (1.9).

Notice that Badia-Duong never made use of the complex plane wave.  Instead, they used constant, first- and second-order 
harmonic polynomials and separated variable harmonic functions.

We think that there should be several applications of our method to the inverse source problems
for other equations and maybe the inverse conductivity problem with a single set of
Cauchy data (see \cite{KS} and the references therein).  In fact, in the final section of this paper we
present an example of one of the applications in such a direction.

\section{Proof of Theorem 1.1}

\proclaim{\noindent Proposition 2.1.}
Let $v$ be the weak solution of (1.1) with $F=0$.  Then
$$\displaystyle
<\frac{\partial u}{\partial\nu}\vert_{\partial\Omega}, v\vert_{\partial\Omega}>
-<\frac{\partial v}{\partial\nu}\vert_{\partial\Omega}, u\vert_{\partial\Omega}>
=\tilde{F}(v).
\tag {2.1}
$$

\endproclaim

This is justfrom the definition of the Neumann derivatives of $u$ and $v$.

Now we substitute a special solution $v$ of $\Delta v+k^2v=0$ in $\Omega$ into (2.1):
$$\displaystyle
v=\exp\{x\cdot(\tau\omega+i\sqrt{\tau^2+k^2}\,\omega^{\perp})-\tau t\}.
$$
Here $\tau>0$ and $t\in\Bbb R$ are parameters and $\omega,\omega\in S^1$ satisfy $\omega\cdot\omega^{\perp}=0$.
Using $v$ indicated above, we have the following definition.

{\bf\noindent Definition (indicator function).}

(i) Let $u$ be the weak solution of (1.1) with $F=\rho\chi_D$.  Define
$$\displaystyle
I_{\omega}(\tau,t)
=
<\frac{\partial u}{\partial\nu}\vert_{\partial\Omega}, v\vert_{\partial\Omega}>
-<\frac{\partial v}{\partial\nu}\vert_{\partial\Omega}, u\vert_{\partial\Omega}>.
$$

(ii) Let $u'$ be the weak solution of (1.1) with $F=\nabla\cdot\{\rho\chi_D\}$.  Define
$$\displaystyle
J_{\omega}(\tau,t)
=
<\frac{\partial u'}{\partial\nu}\vert_{\partial\Omega}, v\vert_{\partial\Omega}>
-<\frac{\partial v}{\partial\nu}\vert_{\partial\Omega}, u'\vert_{\partial\Omega}>.
$$

These indicator functions satisfy the following.

\proclaim{\noindent Proposition 2.2.}
$$\displaystyle
I_{\omega}(\tau,t)=\int_D\rho(x)v(x)\,dx
\tag {2.2}
$$
$$\displaystyle
J_{\omega}(\tau,t)=-\int_D\rho(x)\mbox{\boldmath $a$}\cdot\nabla v(x)\,dx
\tag {2.3}
$$
$$\begin{array}{ll}
\displaystyle
I_{\omega}(\tau,t)=\exp(\tau(s-t))\,I_{\omega}(\tau,s)
& t,s\in\Bbb R
\end{array}
\tag {2.4}
$$
$$
\displaystyle
J_{\omega}(\tau,t)=
-\mbox{\boldmath $a$}\cdot\{\tau\omega+i\sqrt{\tau^2+k^2}\,\omega^{\perp}\}\,I_{\omega}(\tau,t)
\tag {2.5}
$$

\endproclaim

{\it\noindent Proof.}
Recall that for any $\varphi\in H^1_0(\Omega)$
$$\displaystyle
F(\varphi)
=\int_D\rho(x)\varphi(x)\,dx
$$
if $F=\rho(x)\chi_D(x)$;

$$\displaystyle
F(\varphi)
=-\int_D\rho(x)\mbox{\boldmath $a$}\cdot\nabla\varphi(x)\,dx
$$
if $F=\nabla\cdot\{\rho(x)\chi_D(x)\mbox{\boldmath $a$}\}$.

A combination of these facts and (2.1) yields (2.2) and (2.3).
Now (2.4) and (2.5) are trivial.

\noindent
$\Box$

\proclaim{\noindent Lemma 2.1.}
Assume that $\omega$ is regular with respect to $D$ and denote by $p$ the only one point
on the set $\{x\in\Bbb R^2\,\vert\,x\cdot\omega=h_D(\omega)\}\cap\partial D$.
Then there exists a non-zero complex number $A$ such that
$$\displaystyle
I_{\omega}(\tau,h_D(\omega))
=\frac{1}{\tau^2}\tilde{\rho}(p)\exp\{i\sqrt{\tau^2+k^2}\,p\cdot\omega^{\perp}\}\,A+
O\left(\frac{1}{\tau^{2+\alpha}}\right)
\tag {2.6}
$$
as $\tau\longrightarrow\infty$.

\endproclaim

{\it\noindent Proof.}  Take $\eta>0$ and $\tilde{\rho}$ such that (1.7) and (1.8) hold.  Since $p$ is a vertex of $D$ 
and $\omega$ is regular, there exists $\delta>0$ and two linear functions $f(s)$ and $g(s)$ on $[0,\,\delta]$
such that
$$\displaystyle
\begin{array}{l}
\displaystyle
\{x=(h_D(\omega)-s)\omega+y\omega^{\perp}\,\vert\,
g(s)\le y\le f(s),\,0\le s\le\delta\}\\
\\
\displaystyle
=\{x\in\Bbb R^2\,\vert\,h_D(\omega)-\delta\le x\cdot\omega\le h_D(\omega)\}\cap\overline{D}
\subset B(p,\eta)\cap\overline{D}
\end{array}
$$
$$\begin{array}{lll}
\displaystyle
f(0)=g(0)=p\cdot\omega^{\perp}
&
\displaystyle
f'(0)\ge 0,\,g'(0)\le 0
&
\displaystyle
f'(0)-g'(0)>0.
\end{array}
\tag {2.7}
$$
Set
$$\begin{array}{l}
\displaystyle
\tilde{\rho}(s,y)=\tilde{\rho}((h_D(\omega)-s)\omega+y\omega^{\perp})
\,\,\,\,\,\,
I=\int_0^{\delta}e^{-\tau s}
\left(\int_{g(s)}^{f(s)}
e^{i\sqrt{\tau^2+k^2}\,y}\,dy\right)\,ds\\
\\
\displaystyle
II
=\int_0^{\delta}e^{-\tau s}
\left(\int_{g(s)}^{f(s)}
(\tilde{\rho}(s,y)-\tilde{\rho}(0,p\cdot\omega^{\perp}))e^{i\sqrt{\tau^2+k^2}\,y}\,dy\right)\,ds.
\end{array}
$$
Then
$$\begin{array}{ll}
\displaystyle
I_{\omega}(\tau,h_D(\omega))
&
\displaystyle
=\int_0^{\delta}e^{-\tau s}
\left(\int_{g(s)}^{f(s)}
\tilde{\rho}(s,y)e^{i\sqrt{\tau^2+k^2}\,y}\,dy\right)\,ds+O(e^{-\tau\delta})\\
\\
\displaystyle
&
\displaystyle
=\tilde{\rho}(0,p\cdot\omega^{\perp})\,I+II+O(e^{-\tau\delta}).
\end{array}
\tag {2.8}
$$
We get
$$\begin{array}{ll}
\displaystyle
\vert II\vert
&
\displaystyle
\le C\int_0^{\delta}e^{-\tau s}
\left(\int_{g(s)}^{f(s)}
(\vert s\vert^2+\vert y-p\cdot\omega^{\perp}\vert^2)^{\alpha/2}\,dy\right)\,ds
\\
\\
\displaystyle
&
\displaystyle
\le
C'
\int_0^{\delta}
e^{-\tau s}
\left(\int_{g(s)}^{f(s)}
(\vert s\vert^{\alpha}+\vert y-p\cdot\omega^{\perp}\vert^{\alpha})\,dy\right)\,ds.
\end{array}
$$
Since there exists a constant $C>0$ such that, for any $s\in[0,\,\delta]$,
$$\begin{array}{lll}
\displaystyle
0\le f(s)-g(s)\le Cs &
\displaystyle
0\le f(s)-p\cdot\omega^{\perp}\le Cs
&
\displaystyle
0\le p\cdot\omega^{\perp}-g(s)
\le Cs
\end{array}
$$
one gets
$$\displaystyle
\int_{g(s)}^{f(s)}(\vert s\vert^{\alpha}+\vert y-p\cdot\omega^{\perp}\vert^{\alpha})\,dy
\le Cs^{\alpha+1}.
$$
This thus yields
$$\displaystyle
\vert II\vert
\le
C''
\int_0^{\delta}e^{-\tau s}s^{\alpha+1}\,ds=
O\left(\frac{1}{\tau^{2+\alpha}}\right).
\tag {2.9}
$$
On the other hand, integration by parts yields
$$\displaystyle
I=\frac{1}{i\sqrt{\tau^2+k^2}}\,
\int_0^{\delta}\{
e^{-\tau s+i\sqrt{\tau^2+k^2}\,f(s)}
-e^{-\tau s+i\sqrt{\tau^2+k^2}\,g(s)}\}\,ds.
$$
We note that $f'(s)=f'(0)$ and $g'(s)=g'(0)$ since both $f(s)$ and $g(s)$ are linear functions.
Thus
$$\displaystyle
\int_0^{\delta}
e^{-\tau s+i\sqrt{\tau^2+k^2}\,f(s)}\,ds
=\int_0^{\delta}
\frac{\displaystyle (d/ds)e^{-\tau s+i\sqrt{\tau^2+k^2}\,f(s)}}
{\displaystyle -\tau+i\sqrt{\tau^2+k^2}\,f'(0)}\,ds
=\frac{\displaystyle
e^{i\sqrt{\tau^2+k^2}\,p\cdot\omega^{\perp}}}
{\displaystyle
\tau-i\sqrt{\tau^2+k^2}\,f'(0)}
+O\left(\frac{1}{\tau e^{\tau\delta}}\right).
\tag {2.10}
$$
Needless to say, we have the estimate for $g(s)$ similar to (2.10).
A combination of these estimates, (2.8) and (2.9) gives the desired conclusion.  Notice that $A$ is given by the formula
$$\displaystyle
A=\frac{\displaystyle f'(0)-g'(0)}
{\displaystyle
(1-if'(0))(1-ig'(0))}.
\tag {2.11}
$$

\noindent
$\Box$

Theorem 1.1 is just a consequence of the following proposition.

\proclaim{\noindent Proposition 2.3 (reconstruction formula of $h_D(\omega)$).}
Assume that $\omega$ is regular for $D$.  Then for any $t\in\Bbb R$
$$\displaystyle
h_D(\omega)-t
=\lim_{\tau\longrightarrow\infty}
\frac{\displaystyle
\log\vert I_{\omega}(\tau,t)\vert}{\tau}
=
\lim_{\tau\longrightarrow\infty}
\frac{\displaystyle
\log\vert J_{\omega}(\tau,t)\vert}{\tau}.
\tag {2.12}
$$
\endproclaim

{\it\noindent Proof.}
From (2.4) we have
$$\displaystyle
I_{\omega}(\tau,t)
=\exp\{\tau(h_D(\omega)-t)\}I_{\omega}(\tau, h_D(\omega)).
\tag {2.13}
$$
A combination of (2.5), (2.6) and (2.13) gives (2.12).

\noindent
$\Box$

We point out a remarkable fact which is also derived from (2.13) and (2.6).
Assume that $F=\rho(x)\chi_D(x)$, $u\in H^2(\Omega)$ and that $\partial\Omega$ is $C^2$.
We denote by $d\sigma$ the surface measure on $\partial\Omega$.  Then
$$\displaystyle
I_{\omega}(\tau,t)
=\int_{\partial\Omega}
\left(
\frac{\partial u}{\partial\nu}\,v-\frac{\partial v}{\partial\nu}\,u\right)\,d\sigma.
\tag {2.14}
$$
Note that this is not a formal expression.
Give an arbitrary $\epsilon>0$.  Then the Cauchy data of $v$ on the set
$$\displaystyle
\{x\in\Bbb R^2\,\vert\,x\cdot\omega<t-\epsilon\}\cap\partial\Omega
$$
decays exponentially as $\tau\longrightarrow\infty$.  This thus yields the following.

\proclaim{\noindent Proposition 2.4.}
Assume that $\omega$ is regular for $D$.
Then $t>h_D(\omega)$ holds if and only if
$$\displaystyle
\lim_{\tau\longrightarrow\infty}
\tau^2
\int_{\{x\cdot\omega\ge t-\epsilon\}\cap\partial\Omega}
\left(
\frac{\partial u}{\partial\nu}\,v-\frac{\partial v}{\partial\nu}\,u\right)\,d\sigma=0.
$$

\endproclaim

{\it\noindent Proof.}
From (2.13) and (2.14) one gets
$$\begin{array}{ll}
\displaystyle
\tau^2
\int_{\{x\cdot\omega\ge t-\epsilon\}\cap\partial\Omega}
\left(
\frac{\partial u}{\partial\nu}\,v-\frac{\partial v}{\partial\nu}\,u\right)\,d\sigma
&
\displaystyle
=\exp\{\tau(h_D(\omega)-t)\}\tau^2 I_{\omega}(\tau,h_D(\omega))\\
\\
\displaystyle
&
\displaystyle
\,\,\,
-\tau^2
\int_{\{x\cdot\omega\le t-\epsilon\}\cap\partial\Omega}
\left(
\frac{\partial u}{\partial\nu}\,v-\frac{\partial v}{\partial\nu}\,u\right)\,d\sigma.
\end{array}
$$
Now from (2.6) we get the desired conclusions.

\noindent
$\Box$

Proposition 2.4 means that to decide whether $t>h_D(\omega)$ we need only the Cauchy data
of $u$ on the set
$$\displaystyle
\{x\in\Bbb R^2\,\vert\,x\cdot\omega\ge t-\epsilon\}\cap\partial\Omega.
$$

\section{Remark}

The key point for the proof of Theorem 1.1 is to prove the exact algebraic decay of
$$\displaystyle
I_{\omega}(\tau,h_D(\omega))
=\int_D\rho(x)
\exp\{x\cdot(\tau\omega+i\sqrt{\tau^2+k^2}\,\omega^{\perp})-\tau h_D(\omega)\}\,dx
$$
as $\tau\longrightarrow\infty$.
When $\partial D$ is smooth, it is not true.  In this section, we describe an example
for such a phenomenon.  Let $B(p,\epsilon)$ be the open disk with radius $\epsilon$ centred at $p$.
Then we have the following.

\proclaim{\noindent Proposition 3.1.}
For any $t\in\Bbb R$
$$\begin{array}{l}
\displaystyle
\int_{B(p,\epsilon)}
\exp\{x\cdot(\tau\omega+i\sqrt{\tau^2+k^2}\,\omega^{\perp})-\tau t\}\,dx\\
\\
\displaystyle
=2\pi\epsilon^2
\int_0^1
r J_0(k\epsilon r)\,dr
\exp\{p\cdot(\tau\omega+i\sqrt{\tau^2+k^2}\,\omega^{\perp})-\tau t\}
\end{array}
\tag {3.1}
$$
where $J_0$ is the Bessel function of order zero.

\endproclaim

{\it\noindent Proof.}
This is nothing but a consequence of the mean value theorem for the Helmholtz equation
(see, e.g, \cite{CH}) but for the reader's convenience we present a direct proof.
Change of variables yields
$$\begin{array}{l}
\displaystyle
\int_{B(p,\epsilon)}
\exp\{x\cdot(\tau\omega+i\sqrt{\tau^2+k^2}\,\omega^{\perp})-\tau t\}\,dx
=\epsilon^2\exp\{p\cdot(\tau\omega+i\sqrt{\tau^2+k^2}\,\omega^{\perp})-\tau t\}
\\
\\
\displaystyle
\times
\int_{B(0,1)}
\exp\{\epsilon\xi\cdot(\tau\omega+i\sqrt{\tau^2+k^2}\,\omega^{\perp})\}\,d\xi
\end{array}
\tag {3.2}
$$
$$\begin{array}{l}
\displaystyle
\int_{B(0,1)}
\exp\{\epsilon\xi\cdot(\tau\omega+i\sqrt{\tau^2+k^2}\,\omega^{\perp})\}\,d\xi
\\
\\
\displaystyle
=
\int_0^1r
\left(
\int_{S^1}
\exp\{\epsilon r\vartheta\cdot(\tau\omega+i\sqrt{\tau^2+k^2}\,\omega^{\perp})\}\,d\vartheta\right)\,dr
\end{array}
\tag {3.3}
$$
and
$$\begin{array}{l}
\displaystyle
\int_{S^1}
\exp\{\epsilon r\vartheta\cdot(\tau\omega+i\sqrt{\tau^2+k^2}\,\omega^{\perp})\}\,d\vartheta
=\int_0^{2\pi}\exp\{\epsilon r(\tau\cos\theta+i\sqrt{\tau^2+k^2}\,\sin\theta)\}\,d\theta\\
\\
\displaystyle
=\int_{\vert z\vert=1}
\exp\left\{
\epsilon r\left[
\frac{\tau}{2}
\left(z+\frac{1}{z}\right)
+
\frac{\sqrt{\tau^2+k^2}}{2}
\left(z-\frac{1}{z}\right)\right]
\right\}
\,\frac{dz}{iz}.
\end{array}
\tag {3.4}
$$
By elementary calculation, we know that the residue of
$$\displaystyle
\exp\left\{
\epsilon r\left[
\frac{\tau}{2}
\left(z+\frac{1}{z}\right)
+
\frac{\sqrt{\tau^2+k^2}}{2}
\left(z-\frac{1}{z}\right)\right]
\right\}
/z
$$
at $z=0$ is
$$\displaystyle
\sum_{m=0}^{\infty}
\frac{(-1)^m}{(m!)^2}
\left(\frac{k^2\epsilon^2 r^2}{4}\right)^m
=J_0(k\epsilon r).
$$
Combining this with (3.2) to (3.4), we get (3.1).

\noindent
$\Box$

Now assume that $D=B(p,\epsilon)$, $\rho(x)=1$ and $\epsilon$ is small eneough.
From (3.1) we conclude that $I_{\omega}(\tau,h_D(\omega))$ decays exponentially as $\tau\longrightarrow\infty$.
More precisely,
$$\displaystyle
\left\{
t\in\Bbb R\,\vert\,
\lim_{\tau\longrightarrow\infty}I_{\omega}(\tau,t)=0\right\}
=]p\cdot\omega,\,\infty[.
$$
Note that $h_D(\omega)=p\cdot\omega+\epsilon$.  When $D$ is a general domain with smooth boundary, what can one extract from $I_{\omega}(\tau,t)$?

\section{Application to an inverse boundary value problem}

We consider the strong solution $u\in H^2(\Omega)$ of the equation
$$\begin{array}{ll}
\displaystyle
\Delta u+k^2\,n(x)u=0
& \text{in $\Omega$.}
\end{array}
\tag {4.1}
$$
Here $k>0$, $n(x)=1+\chi_D(x)\rho(x)$ and we assume that $(D,\rho)$ satisfies (1.4), (1.5), (1.7) and (1.8) of Section 1 and that
$$\begin{array}{ll}
\displaystyle
\text{$\partial\Omega$ is $C^2$}
&
\displaystyle
\text{$\rho\in L^{\infty}(D)$.}
\end{array}
$$
The problem is to reconstruct $(D,\rho)$ from the Cauchy data $(u\vert_{\partial\Omega}, (\partial/\partial\nu)u\vert_{\partial\Omega})$.
It is closely related to the inverse scattering problem of a time-harmonic plane wave by an absorbing inhomogeneous medium.
In this section we prove that one can reconstruct the convex hull of $D$ from such data.

More precisely, let $\omega\in S^1$ be regular with respect to $D$ and denote by $p=p(\omega)$ the only one point
in the set $\{x\in\Bbb R^2\,\vert\, x\cdot\omega=h_D(\omega)\}\cap\partial D$.
By the Sobolev imbedding theorem, one may assume that for any $\beta\in\,]0,\,1[$
$$\displaystyle
u\in C^{0,\,\beta}(\overline\Omega)
\tag {4.2}
$$
and thus one can consider the value of $u$ at $p$.  The result is as follows.

\proclaim{\noindent Theorem 4.1.}
Assume that
$$\displaystyle
u(p)\not=0.
\tag {4.3}
$$
Then one can calculate $h_D(\omega)$ from the Cauchy data of $u$ on $\partial\Omega$.

\endproclaim

Also from the uniqueness of $h_D(\omega)$ we automatically obtain the following.

\proclaim{\noindent Corollary 3.}  Assume that $(D_1,\rho_1)$ and $(D_2,\rho_2)$ are two pairs satisfying (1.4), (1.5), (1.7) and (1.8).
Let $u_j$ be the strong solution of (4.1) with $(D,\rho)=(D_j,\rho_j)$.
Assume that both $u_1$ and $u_2$ never vanish everywhere in $\Omega$.
Then
$$\displaystyle
\left(u_1\vert_{\partial\Omega},\frac{\partial u_1}{\partial\nu}\vert_{\partial\Omega}\right)
=\left(u_2\vert_{\partial\Omega},\frac{\partial u_2}{\partial\nu}\vert_{\partial\Omega}\right)
$$
implies
$$\displaystyle
[D_1]=[D_2].
$$

\endproclaim

The non-vanishing of $u$ everywhere in $\Omega$ can be realized in several situations.
For example, when one consideres $u$ the restriction to $\Omega$ of the total field of the scattering
solution by the incident plane wave $e^{ikx\cdot d}$, $d\in S^1$ of the equation $\Delta u+k^2\,n(x)u=0$ in $\Bbb R^2$,
$u$ never vanishes for $k$ sufficiently small.

{\it\noindent Proof of Theorem 4.1.}
Let $v$ be the strong solution of the Helmholtz equation $\Delta v+k^2 v=0$ in $\Omega$.
Integration by parts yields
$$\displaystyle
\int_{\partial\Omega}
\left(\frac{\partial u}{\partial\nu}\,v-u\,\frac{\partial v}{\partial\nu}\right)\,d\sigma
=-k^2\int_D\rho(x)u(x)v(x)\,dx.
\tag {4.4}
$$
Set $\rho'(x)=\rho(x)u(x)\,(\in L^2(D))$.
From (1.7), (1.8), (4.2) and (4.3) we know that $\rho'$ satisfies (1.7) and (1.8).  Now take
$$\displaystyle
v=\exp\{x\cdot(\tau\omega+i\sqrt{\tau^2+k^2}\,\omega^{\perp})-\tau t\}
\,\,\,\,\,\,\tau>0,\,t\in\Bbb R
$$
and define
$$\displaystyle
I_{\omega}(\tau,t)
=\int_{\partial\Omega}
\left(\frac{\partial u}{\partial\nu}\,v-u\,\frac{\partial v}{\partial\nu}\right)\,d\sigma.
$$
On account of (4.3) a combination of (2.6) and (4.4) as in Proposition 2.3 gives
$$\displaystyle
h_D(\omega)-t
=\lim_{\tau\longrightarrow\infty}
\frac{\log\vert I_{\omega}(\tau,t)\vert}{\tau}
\,\,\,\,\,\,\forall\,t\in\Bbb R
\tag {4.5}
$$
and also
$$\displaystyle
\lim_{\tau\longrightarrow\infty}\tau^2\int_{\{x\cdot\omega\ge t-\epsilon\}\cap\partial\Omega}
\left(\frac{\partial u}{\partial\nu}\,v-\frac{\partial v}{\partial\nu}\,u\right)\,d\sigma
=0
\tag {4.6}
$$
if and only if $t>h_D(\omega)$.

\noindent
$\Box$

$$\quad$$

\centerline{{\bf Acknowledgments}}

The author thanks the referees for several suggestions for the improvement of the manuscript.

$$\quad$$

\end{document}